\documentclass[12pt,reqno]{amsart} 
\usepackage{amsthm,amsfonts,amssymb,amscd}

\newtheorem{theorem}[subsection]{Theorem}
\newtheorem{proposition}[subsection]{Proposition}

\newtheorem{lemma}[subsection]{Lemma}
\newtheorem{corollary}[subsection]{Corollary}

\theoremstyle{definition}
\newtheorem{definition}[subsection]{Definition}

\newtheorem{proposition-definition}[subsection]{Proposition-Definition}

\theoremstyle{remark}

\newcommand{\dual}{{\scriptscriptstyle\vee}}
\newcommand{\mt}[1]{\operatorname{#1}}

\newcommand{\rk}{\operatorname{rk}\nolimits}
\newcommand{\Blowup}{\operatorname{Blowup}\nolimits}

\newcommand{\id}{\operatorname{id}\nolimits}

\newcommand{\Ext}{\operatorname{Ext}\nolimits}

\newcommand{\Cl}{\operatorname{Cl}}
\newcommand{\Pl}{\operatorname{Pl}}

\newcommand{\CC}{{\mathbb C}}

\newcommand{\ZZ}{{\mathbb Z}}

\newcommand{\PP}{{\mathbb P}}
\newcommand{\PPPP}{{{\mathbb P}^{14}}}

\newcommand{\FF}{{\mathbb F}}

\newcommand{\OOO}{{\mathcal O}}

\newcommand{\III}{{\mathcal I}}

\newcommand{\EEE}{{\mathcal E}}
\newcommand{\HHH}{{\mathcal H}}

\newcommand{\GGG}{{\mathcal G}}
\newcommand{\BBB}{{\mathcal B}}
\newcommand{\CCC}{{\mathcal C}}
\newcommand{\KKK}{{\mathcal K}}
\newcommand{\NNN}{{\mathcal N}}

\newcommand{\QQQ}{{\mathcal Q}}
\newcommand{\SSS}{{\mathcal S}}

\newcommand{\EXT}{{\mathcal Ext}}

\newcommand{\LR}{\Longleftrightarrow}
\newcommand\inner{
\begin{picture}(9,8)
\put(1,0){\line(1,0){6}}
\put(1,0){\line(0,1){7}}
\end{picture}
}
\newcommand\surj{{\twoheadrightarrow}}

\newcommand\alp{\alpha}

\newcommand\eps{\epsilon}

\newcommand\la{\lambda}
\newcommand\La{\Lambda}
\newcommand\G{\Gamma}

\newcommand\si{\sigma}

\newcommand\lra{{\longrightarrow}}

\newcommand\Hilb{{\operatorname{Hilb}\nolimits^{5n}_X}}

\renewcommand\square{\frame{\phantom{{\large x}}}}
\renewcommand\empty{\varnothing}

\author{A. Iliev}

\address{A. I.: Inst. of Math.,
Bulgarian Acad. of Sci.,
Acad. G. Bonchev Str., 8,\ 
1113 Sofia, Bulgaria}
\email{ailiev@math.bas.bg}

\author{D. Markushevich}

\address{D. M.: Math\'ematiques - b\^{a}t. M2, Universit\'e Lille 1,
F-59655 Villeneuve d'Ascq Cedex, France}
\email{markushe@gat.univ-lille1.fr}


\title{The Abel--Jacobi map
for \\ a cubic threefold and periods of Fano \\
threefolds of degree 14}

\begin{document}

\begin{abstract}
The Abel--Jacobi maps of the families of elliptic quintics and
rational quartics lying
on a smooth cubic threefold are studied. It is proved that
their generic fiber is the 5-dimensional
projective space for quintics, and a smooth 3-dimensional
variety birational to the cubic itself for quartics.
The paper is a continuation of the recent work of
Markushevich--Tikhomirov, who showed that the first Abel--Jacobi
map factors through the moduli component of stable rank 2 vector
bundles on the cubic threefold with Chern numbers $c_1=0,
c_2=2$ obtained by Serre's construction from elliptic quintics,
and that the factorizing map from the moduli space to the intermediate Jacobian
is \'etale. The above result implies that the degree of the \'etale map
is 1, hence the moduli component of vector bundles
is birational to the intermediate Jacobian. As an applicaton,
it is shown that the generic fiber of the period map of Fano varieties
of degree 14 is birational to the intermediate Jacobian of the
associated cubic threefold.

\end{abstract}

\maketitle

\section*{Introduction}

Clemens and Griffiths studied in \cite{CG} the Abel--Jacobi
map of the family of lines on a cubic threefold $X$. They represented
its intermediate Jacobian $J^2(X)$ as the Albanese variety 
$\mt{Alb} F(X)$ of the Fano surface $F(X)$ parametrizing
lines on $X$ and described its theta divisor. From this description, they deduced
the Torelli Theorem and the non-rationality of $X$. 
Similar results were obtained by Tyurin \cite{Tyu} and Beauville
\cite{B}. 

One can easily understand the structure of the Abel--Jacobi maps
of some other familes of curves of low degree on $X$ (conics,
cubics or elliptic quartics), in reducing
the problem to the results of Clemens--Griffiths and Tyurin.
The first non trivial cases are those of rational normal quartics
and of elliptic normal quintics. We determine the fibers
of the Abel--Jacobi maps of these families of curves, in continuing
the work started in \cite{MT}.

Our result on elliptic quintics
implies that the moduli space of
instanton vector bundles of charge 2 on $X$ has a component,
birational to $J^2(X)$. We conjecture that the moduli space
is irreducible, but the problem of irreducibility stays beyond the scope
of the present article. As far as we know, this is the first example
of a moduli space of vector bundles which is birational to an abelian
variety, different from the Picard or Albanese variety of the base.
The situation is also quite different from the known cases where
the base is $\PP^3$ or the 3-dimensional quadric. In these cases,
the instanton moduli space is irreducible and rational
at least for small charges,
see \cite{Barth}, \cite{ES}, \cite{H}, \cite{LP}, \cite{OSz}.
Remark, that for the cubic $X$, two is the smallest possible charge,
but the moduli space is non-rational.
There are no papers on the geometry of particular
moduli spaces of vector bundles
for other 3-dimensional Fano varieties (for some constructions
of vector bundles on such varieties, see \cite{G1}, \cite{G2}, \cite{BM},
\cite{BM2},\cite{SW}, \cite{A-C}).

The authors of \cite{MT} proved that
the Abel--Jacobi map $\Phi $ of the family of elliptic quintics lying
on a general cubic threefold $X$ factors through a 5-dimensional moduli 
component $M_X$ of stable rank 2 vector
bundles $\EEE$ on $X$ with Chern numbers $c_1=0,
c_2=2$. The factorizing map $\phi$ sends an elliptic
quintic $C\subset X$ to
the vector bundle $\EEE$ obtained by Serre's construction from $C$ (see Sect. 2). 
The fiber $\phi^{-1}([\EEE ])$ is a 5-dimensional projective space
in the Hilbert scheme $\Hilb$,
and the map $\Psi$ from the moduli space to the intermediate Jacobian $J^2(X)$,
defined by $\Phi =\Psi\circ\phi$,
is \'etale on the open set representing (smooth) elliptic quintics which are
not contained in a hyperplane (Theorem \ref{mainMT}).

We improve the result of \cite{MT} in showing that the degree of the
above \'etale map is 1. Hence $M_X$ is birational to $J^2(X)$ and
the generic fiber of $\Phi$ is just
one copy of $\PP^5$ (see Theorem \ref{degree1} and Corollary \ref{openinJ}).
The behavior of the Abel--Jacobi map of
elliptic quintics is thus quite similar to that of the Abel--Jacobi
map of divisors on a curve, where all the fibers are projective spaces.
But we prove that the situation is very different in the case of
rational normal quartics, where the fiber of the Abel--Jacobi map
is a {\em non-rational} 3-dimensional variety: it is birationally equivalent to
the cubic $X$ itself
(Theorem \ref{fiberquart}). 

The first new ingredient of our proofs, comparing to \cite{MT},
is another interpretation of the vector bundles $\EEE $ from $M_X$.
We represent the cubic $X$ as a linear section of the Pfaffian
cubic in $\PP^{14}$, parametrizing $6\times 6$ matrices $M$ of rank 4, and
realize $\EEE^\dual (-1)$
as the restriction of the kernel bundle $M\mapsto\ker M\subset\CC^6$
(Theorem \ref{KKK}).
The kernel bundle has been investigated by A. Adler in his Appendix to
\cite{AR}. We prove that it embeds $X$ into the Grassmannian $G=G(2,6)$,
and the quintics $C\in\phi^{-1}([\EEE ])$ become the sections of $X$
by the Schubert varieties $\si_{11}(L)$ for all hyperplanes $L\subset \CC^6$.
We deduce that for any line $l\subset X$, each fiber of $\phi$ contains 
precisely one pencil $\PP^1$ of reducible curves of the form $C'+l$
(Lemma \ref{Hliz}).
Next we use the techniques of Hartshorne--Hirschowitz \cite{HH} for smoothing
the curves of the type ``a rational normal quartic plus one of its chords in $X$"
(see Sect. 4)
to show that there is a 3-dimensional family
of such curves in a generic fiber of $\phi$  and that the above pencil $\PP^1$ 
for a generic $l$
contains curves $C'+l$ of this type (Lemma \ref{dim3}, Corollary \ref{rnqpc}).

The other main ingredient is the parametrization of $J^2(X)$
by minimal sections of the 2-dimensional conic bundles of the form
$Y(C^2)=\pi_l^{-1}(C^2)$, where $\pi_l:\mt{Blowup}_l(X)\lra\PP^2$
is the conic bundle obtained by projecting $X$ from a fixed line $l$,
and $C^2$ is a generic conic in $\PP^2$ (see Sect. 3).
The standard Wirtinger approach \cite{B}
parametrizes $J^2(X)$ by reducible curves which are sums of components
of reducible fibers of $\pi_l$. Our approach, developed
in \cite{Ili} in a more general form, replaces the degree 10
sums of components of the reducible fibers
of the surfaces  $Y(C^2)$ by the irreducible curves which are sections
of the projection $Y(C^2)\lra C^2$ with a certain minimality condition.
This gives a parametrization of $J^2(X)$ by a family of rational curves,
each one of which is projected isomorphically onto some conic in $\PP^2$.
It turns out, that these rational curves are normal quartics meeting
$l$ at two points. They form a {\it unique} pencil $\PP^1$
in each fiber of the Abel--Jacobi map of rational normal quartics.
Combining this with the above, we conclude that the curves of type
$C'+l$ form a unique pencil in each fiber of $\Phi$, hence the
fiber is one copy of $\PP^5$.

In conclusion, we provide a description of the moduli space
of Fano varieties $V_{14}$ as a birationally fibered space
over the moduli space of cubic 3-folds with the intermediate
Jacobian as a fiber (see Theorem \ref{periods}). The interplay
between cubics and varieites $V_{14}$ is exploited several times
in the paper. We use the Fano--Iskovskikh birationality between 
$X$ and $V_{14}$ to prove Theorem \ref{KKK} on kernel bundles, 
and the Tregub--Takeuchi
one (see Sect. 1) to study the fiber of the Abel--Jacobi map of the family of
rational quartics (Theorem \ref{fiberquart}) and the
relation of this family to that of normal elliptic quintics
(Proposition \ref{claim}).

\medskip

{\em Acknowledgements}. The second author acknowledges with
pleasure the hospitality
of the MPIM at Bonn, where
he completed the work on the paper.

\section{Birational isomorphisms between $V_3$ and $V_{14}$}

There are two constructions of birational isomorphisms between
a nonsingular cubic threefold $V_3\in\PP^4$ and the Fano variety $V_{14}$
of degree 14 and of index 1,
which is a nonsingular section of the Grassmannian $G(2,6)\in\PP^{14}$
by a linear subspace of codimension 5. The first one is that of Fano--Iskovskikh,
and it gives a birational isomorphism whose indeterminacy locus in both 
varieties is an elliptic curve together with some 25 lines;
the other is due to Tregub--Takeuchi,
and its indeterminacy locus is a rational quartic plus 16 lines
on the side of $V_3$,
and 16 conics passing through one point on the side of $V_{14}$.
We will sketch both of them.

\begin{theorem}[Fano--Iskovskikh]\label{FI}
Let $X=V_3$ be a smooth cubic threefold. Then $X$ contains
a smooth projectively normal elliptic quintic curve. Let $C$ be such a curve.
Then $C$ has exactly 25 bisecant
lines $l_i \subset X$, $i = 1,...,25$, and
there is a unique effective divisor
$M \in \mid {\mathcal O}_X(5 - 3C) \mid$ on $X$,
which is a reduced surface containing the $l_i$. The following assertions
hold:

(i) 
The non-complete linear system
$\mid {\mathcal O}_X(7 - 4C) \mid$ defines a birational map
${\rho}:X \rightarrow V$ where $V=V_{14}$ is a Fano 3-fold of
index 1 and of degree 14.
Moreover $\rho = \sigma \circ \kappa \circ \tau$
where ${\sigma}:X' \rightarrow X$ is the blow-up of $C$,
${\kappa}:X' \rightarrow X^+$ is a flop over the proper
transforms $l_i' \subset X'$ of the $l_i$, $i=1,...,25$,
and ${\tau}:X^+ \rightarrow V$ is a blowdown of
the proper transform $M^+ \subset X^+$ of $M$
onto an elliptic quintic $B \subset V$.
The map $\tau$ sends the transforms $l_i^+ \subset X^+$
of $l_i$ to the 25 secant lines $m_i \subset V$, $i=1,...,25$
of the curve $B$.

(ii)
The inverse map ${\rho}^{-1}$ is defined by the system
$\mid {\mathcal O}_V(3 - 4B) \mid$. The exceptional divisor
$E' = {\sigma}^{-1}(C) \subset X'$ is the proper transform
of the unique effective divisor
$N \in \mid {\mathcal O}_V(2 - 3B) \mid$.
\end{theorem}

For   a proof , see \cite{Isk-1}, \cite{F}, or \cite{IP}, Ch. 4.

\begin{theorem}[Tregub--Takeuchi]\label{TT}
Let $X$ be a smooth cubic threefold. Then $X$ contains 
a rational projectively normal quartic
curve. Let $\G$ be such a curve.
Then $\G $ has exactly 16 bisecant
lines $l_i \subset X$, $i = 1,...,16$,and there is a unique effective divisor
$M \in \mid {\mathcal O}_X(3 - 2\G ) \mid$ on $X$,
which is a reduced surface containing the $l_i$.
The following assertions hold:

(i)
The non-complete linear system
$\mid {\mathcal O}_X(8 - 5\G ) \mid$ defines a birational map
${\chi}:X \rightarrow V$ where $V$ is a Fano 3-fold of
index 1 and of degree 14.
Moreover $\chi = \sigma \circ \kappa \circ \tau$,
where ${\sigma}:X' \rightarrow X$ is the blowup of $\G $,
${\kappa}:X' \rightarrow X^+$ is a flop over the proper
transforms $l_i' \subset X'$ of $l_i$, $i=1,...,16$,
and ${\tau}:X^+ \rightarrow V$ is a blowdown of
the proper transform $M^+ \subset X^+$ of $M$
to a point $P \in V$.
The map $\tau$ sends the transforms $l_i^+ \subset X^+$
of $l_i$ to the 16 conics $q_i \subset V$, $i=1,...,16$
which pass through the point $P$.

(ii)
The inverse map ${\chi}^{-1}$ is defined by the system
$\mid {\mathcal O}_V(2 - 5P) \mid$. 
The exceptional divisor
$E' = {\sigma}^{-1}(\G ) \subset X'$ is the proper transform
of the unique effective divisor
$N \in \mid {\mathcal O}_V(3 - 8P) \mid$.

(iii) For a generic
point $P$ on any nonsingular $V_{14}$, this linear system
defines a birational isomorphism of type ${\chi}^{-1}$.

\end{theorem}

\begin{proof} For (i), (ii), see \cite{Tak}, Theorem 3.1, and \cite{Tre}.
For (iii), see \cite{Tak}, Theorem 2.1, (iv). See also \cite{IP}, Ch. 4.
\end{proof}

\subsection{Geometric description}\label{geom-puts}

We will briefly describe the geometry of the first birational isomorphism
between $V_3$ and $V_{14}$ following \cite{Pu}.

Let $E$ be a 6-dimensional vector space over $\CC$. Fix a basis
$e_0,\ldots ,e_5$
for  $E$, then $e_i\wedge e_j$ for
$0\leq i<j\leq 5$ form a basis for the Pl\"ucker space of 2-spaces in $E$,
or equivalently, 
of lines in
$\PP^5=\PP (E)$.  With Pl\"ucker coordinates $x_{ij}$,
the embedding of the Grassmannian $G=G(2,E)$ in 
$\PP^{14}=\PP 
(\wedge^2E)$ is precisely the locus of rank 2 skew symmetric
$6\times 6$ matrices
$$ M=\left[\begin{array}{cccccc} 0&x_{01}&x_{02}&x_{03}&x_{04}&x_{05}\\
-x_{01}&0&x_{12}&x_{13}&x_{14}&x_{15}\\
-x_{02}&-x_{12}&0&x_{23}&x_{24}&x_{25}\\
-x_{03}&-x_{13}&-x_{23}&0&x_{34}&x_{35}\\
-x_{04}&-x_{14}&-x_{24}&-x_{34}&0&x_{45}\\
-x_{05}&-x_{15}&-x_{25}&-x_{35}&-x_{45}&0\\ \end{array}\right]
.$$

There are two ways to associate to these data a 13-dimensional cubic.
The Pfaffian cubic hypersurface $\Xi \subset\PP^{14}$ is defined as the
zero locus of the $6\times 6$ Pfaffian of this
matrix; it can be identified with the secant variety of $G(2,E)$, or else,
it is the locus where $M$ has rank 4. The other way is to consider
the dual variety $\Xi '=G^\dual \subset\PP^{14\dual }$ of $G$; it is also a cubic
hypersurface, which is nothing other than the secant variety of the
Grassmannian $G'=G(2,E^\dual )\subset\PP (\wedge^2E^\dual )=\PP^{14\dual}$.

As it is classically known,
the generic cubic threefold $X$ can be represented as a section of the
Pfaffian cubic by a linear subspace of codimension 10; see also
a recent proof  in \cite{AR}, Theorem 47.3. There are $\infty^5$ essentially different ways to do this.
Beauville and Donagi \cite{BD}
have used this idea for introducing the symplectic
structure on the Fano 4-fold (parametrizing lines) of a cubic
4-fold. In their case, only special cubics (a divisorial family)
are sections of the Pfaffian cubic, so they introduced the symplectic structure
on the Fano 4-folds of these special cubics, and obtained the existence
of such a structure on the generic one by deformation arguments.

For any hyperplane section $H\cap G$ of $G$, we can define $\rk H$ as the rank
of the antisymmetric matrix $(\alp_{ij})$, where $\sum \alp_{ij}x_{ij}=0$
is the equation of $H$. So, $\rk H$ may take the values 2,4 or 6.
If $\rk H=6$, then $H\cap G$ is nonsingular and
for any $p\in\PP^5=\PP (E)$, there is the unique
hyperplane $L_p\subset\PP^5=\PP (E)$, such that $q\in H\cap G$, $p\in l_q$
$\LR$ $l_q\subset L_p$. Here $l_q$ denotes the line in $\PP^5$
represented by $q\in G$. (This is a way to see that the base of the
family of 3-dimensional planes on the 7-fold $H\cap G$ is $\PP^5$.)

The rank of $H$ is 4 if and only if $H$ is tangent to $G$ at exactly
one point $z$, and in this case, the hyperplane $L_p$ is not defined
for any $p\in l_z$: we have for such $p$ the equivalence
$p\in l_x$ $\LR$ $x\in H$. Following Puts, we call the line
$l_z$ the {\it center} of $H$; it will be denoted $c_H$.

In the third case, when $\rk H=2$,
$H\cap G$ is singular along the whole Grassmannian 
subvariety $G(2,4)=G(2,E_H)$, where $E_H=\ker (\alp_{ij})$
is of dimension 4. We have $x\in H$ $\LR$ $l_x\cap \PP(E_H)\neq\empty$.

This description identifies the dual of $G$ with
$\Xi '=\{ H\mid\rk H\leq 4\} =\{ H\mid P\! f ((\alp_{ij}))=0$,
and its singular locus with $\{E_H\}_{\rk H=2}=G(4,E)$.

Now, associate to any nonsingular 
$V_{14}=G\cap \Lambda$, where $\Lambda =H_1\cap H_2\cap H_3\cap H_4\cap H_5$,
the cubic 3-fold $V_3$ by the following rule:
\begin{equation}\label{V3toV14}
V_{14}=G\cap \Lambda\
\mapsto \ V_3=
\Xi '\cap \Lambda^\dual , 
\end{equation}
where $\Lambda^\dual =<H_1^\dual ,
H_2^\dual ,H_3^\dual ,H_4^\dual ,H_5^\dual >$, $H_i^\dual$ denotes
the orthogonal complement of $H_i$ in $\PP^{14\dual }$,
and the angular brackets the linear span.  One can prove
that $V_3$ is also nonsingular.

According to Fano, the lines $l_x$ represented by points $x\in V_{14}$
sweep out an irreducible quartic hypersurface $W$, which Fano calls
the quartic da Palatin. $W$ coincides with the union of
centers of all $H\in V_3$.
One can see, that $W$ is singular
along the locus of foci $p$ of Schubert pencils of lines on $G$
$$\sigma_{43}(p,h)=\{ x\in G\mid p\in l_x\subset h\} 
$$
which lie entirely in $V_{14}$, where $h$ denotes a plane in $\PP^5$
(depending on $p$). The pencils $\sigma_{43}$ are exactly the lines
on $V_{14}$, so Sing$\, W$ is identified with the base of the
family of lines on $V_{14}$, which is known to be a nonsingular curve
of genus 26 for generic $V_{14}$ (see, e. g. \cite{M} for the study
of the curve of lines on $V_{14}$, and Sections
50, 51 of \cite{AR} for the study of Sing$\, W$ without any connection to
$V_{14}$).

The construction of the birational isomorphism $\eta_L :V_{14}
\dasharrow V_3$ depends on the choice of a hyperplane $L\subset\PP^5$.
Let 
$$
\phi :V_{14}\dasharrow W\cap L, \  x\mapsto L\cap l_x \ ,
\psi :V_3\dasharrow W\cap L, \  H^\dual\mapsto L\cap c_H .
$$
These two maps
are birational, and $\eta_L$ is defined by 
\begin{equation}\label{chi}
\eta_L = \psi^{-1}\circ\phi .
\end{equation}

The locus, on which $\eta_L$ is not an isomorphism, consists
of points where either $\phi$ or $\psi$ is not defined or is
not one-to-one. The indeterminacy locus $B$ of $\phi$ consists of all the points $x$ such that $l_x\subset L$, that
is, $B=G(2,L)\cap H_1\cap\ldots\cap H_5$. For generic $L$, it is
obviously a smooth elliptic quintic curve in $V_{14}$, and it is this
curve that was denoted in Theorem \ref{FI} by the same symbol $B$. The indeterminacy locus of $\psi$ is described
in a similar way.
We summarize the above in the following statement.

\begin{proposition}\label{puts}
Any nonsingular variety $V_{14}$ determines a unique nonsingular cubic $V_3$
by the rule (\ref{V3toV14}). Conversely, a generic cubic $V_3$
can be obtained in this way from $\infty^5$  many varieties $V_{14}$.

For each pair $(V_{14}, V_3)$ related by (\ref{V3toV14}),
there is a family of birational maps $\eta_L :V_{14}
\dasharrow V_3$, defined by (\ref{chi}) and parametrized
by points of the dual projective
space $\PP^{5\dual }$, and the structure of $\eta_L$ for generic $L$ is described by Theorem \ref{FI}.

The smooth elliptic quintic 
curve $B$ (resp. $C$) of Theorem \ref{FI} is the locus of points $x\in V_{14}$
such that $l_x\subset L$ (resp. $H^\dual\in V_3$
such that $c_H\subset L$).
\end{proposition}

\begin{definition}
We will call two varieties $V_3$, $V_{14}$ associated (to each other),
if $V_3$ can be obtained from $V_{14}$ by the construction (\ref{V3toV14}).
\end{definition}

\subsection{Intermediate Jacobians of $V_3$, $V_{14}$}

Both constructions of birational isomorphisms give the isomorphism
of the intermediate Jacobians of generic varieties $V_3$,
$V_{14}$, associated to each other.
This is completely obvious for
the second construction: it gives a birational isomorphism, which
is a composition of blowups and blowdowns with centers in nonsingular rational
curves or points. According to \cite{CG}, a blowup 
$\sigma :\tilde{X}\lra X$ of a threefold $X$ with a
nonsingular center $Z$ can change its intermediate Jacobian only in the
case when $Z$ is a curve of genus $\geq 1$, and in this case
$J^2(\tilde{X})\simeq J^2(X)\times J(Z)$ as principally polarized abelian
varieties, where $J^2$ (resp. $J$) stands for the intermediate Jacobian
of a threefold (resp. for the Jacobian of a curve). Thus, the Tregub--Takeuchi
birational isomorphism does not change the intermediate Jacobian.
Similar argument works for the Fano--Iskovskikh construction. It factors
through blowups and blowdowns with centers in rational curves,
and contains in its factorization exactly one blowup and one blowdown
with nonrational centers, which are elliptic curves. So, we have
$J^2(V_3)\times C\simeq J^2(V_{14})\times B$ for some elliptic curves
$C,B$.  According to Clemens--Griffiths,  $J^2(V_3)$
is irreducible for every nonsingular $V_3$, so we can simplify the
above isomorphism to obtain $J^2(V_3)\simeq J^2(V_{14})$ (and we also
obtain, as a by-product, the isomorphism $C\simeq B$).

\begin{proposition}\label{FI-TT}
Let $V=V_{14}$, $X=V_3$ be a pair of smooth Fano 
varieties related by
either of the two birational isomorphisms of Fano--Iskovskikh or of
Tregub--Takeuchi. Then $J^2(X)\simeq J^2(V)$,  $V,X$ are associated to each other
and related by a birational isomorphism of the other type as well.

\end{proposition}

\begin{proof}
The isomorphism of the intermediate Jacobians was proved in the previous
paragraph.
Let $J^2(V') = J^2(V'') = J$. By Clemens-Griffiths \cite{CG}
or Tyurin \cite{Tyu}, the global Torelli Theorem holds for smooth
3-dimensional cubics, so there exists
the unique cubic threefold $X$ such that
$J^2(X) = J$ as p.p.a.v.
Let $X'$ and $X''$ be the unique cubics associated to $V'$ and $V''$.
Since $J^2(X') = J^2(V') = J = J^2(V'') = J^2(X'')$, then
$X'\simeq X\simeq  X''$.

Let now $V'$ and $V''$ be associated to the same cubic threefold $X$,
and let $J^2(X) = J$.
Then by the above $J^2(V') = J^2(X) = J^2(V'')$.

Let $X$, $V$ be related
by, say, a Tregub--Takeuchi birational isomorphism. By
Proposition \ref{puts}, $V$ contains a smooth elliptiic quintic curve
and admits a birational isomorphism of Fano--Iskovskikh type
with some cubic $X'$. Then, as above, $X\simeq X'$ by Global Torelli,
and $X$, $V$ are associated to each other by the definition of
the Fano--Iskovskikh birational isomorphism. Conversely,
if we start from the hypothesis that $X$, $V$
are related by a Fano--Iskovskikh birational isomorphism, then
the existence of a Tregub--Takeuchi one from $V$ to some
cubic $X'$ is affirmed by Theorem \ref{TT}, (iii). Hence,
again by Global Torelli, $X\simeq X'$ and we are done.

\end{proof}

\section{Abel--Jacobi map and
vector bundles on a cubic threefold}

Let   $X$ be a smooth cubic threefold. The authors of \cite{MT} have 
associated to every normal elliptic quintic curve
$C\subset X$ a stable rank 2 vector bundle $\EEE =
\EEE_C$, unique up to isomorpfism. It is defined
by Serre's construction:
\begin{equation}\label{serre}
0\lra \OOO_X\lra \EEE (1) \lra \III_C(2) \lra 0\; ,
\end{equation}
where $\III_C=\III_{C,X}$ is the ideal sheaf of $C$ in $X$.
Since the class of $C$ modulo algebraic equivalence
is $5l$, where $l$ is the class of a line,
the sequence (\ref{serre}) implies that
$c_1(\EEE )=0, c_2(\EEE )=2l$. One sees immediately from
(\ref{serre}) that $\det\EEE$ is trivial,
and hence  $\EEE$ is self-dual as soon as it is a vector
bundle (that is, $\EEE^\dual \simeq\EEE$).
See \cite[Sect. 2]{MT} for further details on this construction.

Let $\HHH^* \subset\Hilb$ be the open set of the Hilbert scheme 
parametrizing normal elliptic quintic curves in $X$, and
$M\subset M_X(2;0,2)$ the open subset in the moduli space of
vector bundles on $X$ parametrizing those stable
rank 2 vector bundles which arise via Serre's construction
from normal elliptic quintic curves. Let  $\phi^* :\HHH^*
\lra M$ be the natural map. For any reference curve $C_0$ of degree 5
in $X$, let $\Phi^* :\HHH^*\lra J^2(X)$,
$[C]\mapsto [C-C_0]$, be the Abel--Jacobi map.
The following result is proved in \cite{MT}.

\begin{theorem}\label{mainMT}
$\HHH^*$ and $M$ are smooth of dimensions 10 and 5
respectively. They are also irreducible for generic $X$.
There exist a bigger open subset $\HHH\subset\Hilb$ in the
nonsingular locus of $\Hilb$ containing
$\HHH^*$ as a dense subset and extensions of $\phi^*,\Phi^*$ to morphisms
$\phi , \Phi$ respectively, defined on the whole of $\HHH$,
such that the following properties are verified:

(i) $\phi$ is a locally trivial fiber bundle in the \'etale
topology with fiber $\PP^5$. For every $[\EEE ]\in M$, we have
$h^0 (\EEE (1))=6$, and
$\phi^{-1}([\EEE ])\subset \HHH$ is nothing but the $\PP^5$ of
zero loci of all the sections of $\EEE (1)$.

(ii) The fibers of $\Phi$ are finite unions of those of
$\phi$, and the map $\Psi :M\lra J^2(X)$ in the natural factorization
$\Phi =\Psi\circ\phi$ is a quasi-finite \'etale morphism.
\end{theorem}

Now, we will give another interpretation of the vector bundles
$\EEE_C$. Let us represent the cubic $X=V_3$ as a section of the
Pfaffian cubic $\Xi '\subset \PPPP^\dual $ and keep the notation of
\ref{geom-puts}. Let $\KKK$ be the kernel bundle on $X$
whose fiber at $M \in X$ is $\ker H$. Thus $\KKK$ is a
rank 2 vector subbundle of the trivial rank 6 vector bundle
$E_X=E\otimes_\CC \OOO_X$. Let $i :X\lra\PPPP$ be the composition
$\Pl\circ\Cl$, where $\Cl :X\lra G(2,E)$ is the classifying map
of $\KKK\subset E_X$, and $\Pl :G(2,E)\hookrightarrow
\PP (\wedge^2E)=\PPPP$ the Pl\"ucker embedding.

\begin{theorem}\label{KKK}
For any vector bundle $\EEE$ obtained by Serre's construction
starting from a normal elliptic quintic $C\subset X$, there exists
a representation of $X$ as a linear section of $\Xi '$ such that 
$\EEE (1)\simeq\KKK^\dual$ and all
the global sections of $\EEE (1)$ are the images of the constant sections
of $E_X^\dual$ via the natural map $E_X^\dual\lra\KKK^\dual$.
For generic $X,\EEE$, such a representation is unique modulo
the action of $PGL(6)$ and the map $i$ can be identified with
the restriction $v_2|_X$ of the
Veronese embedding $v_2 : \PP^4\lra\PPPP$ of degree 2.
\end{theorem}

\begin{proof}
Let $C\subset X$ be a normal elliptic quintic. By Theorem
\ref{FI}, there exists a $V_{14}=G\cap\Lambda$ together with a birational isomorphism
$X\dasharrow V_{14}$. Proposition \ref{FI-TT} implies that $X$ and $V_{14}$
are associated to each other. By Proposition \ref{puts}, we have $C=
\{ H^\dual\in X\mid c_H\subset L\} =\Cl^{-1}(\si_{11}(L))$, where
$\si_{11}(L)$ denotes the Schubert variety in $G$ parametrizing the
lines $c\subset\PP (E)$ contained in $L$. It is standard that $\si_{11}(L)$
is the scheme of zeros of a section of the dualized universal rank
2 vector bundle $\SSS^\dual$ on $G$. Hence $C$ is the scheme of zeros of
a section of $\KKK^\dual =\Cl^*(\SSS^\dual )$. Hence $\KKK^\dual$
can be obtained by Serre's construction from $C$, and by uniqueness,
$\KKK^\dual \simeq \EEE_C(1)$.

By Lemma 2.1, c) of \cite{MT}, $h^0(\EEE_C(1))=6$, so, to prove the assertion
about global sections, it is enough to show the injectivity of the natural
map $E^\dual =H^0(E_X^\dual )\lra H^0(\KKK^\dual )$. The latter is obvious,
because the quartic da Palatini is not contained in a hyperplane. Thus
we have $E^\dual =H^0(\KKK^\dual )$.

For the identification of $i$ with $v_2|_X$, it is sufficient to
show that $i$ is defined by the sections of $\OOO (2)$ in the image of 
the map $\mt{ev}:\Lambda^2H^0(\EEE (1))\lra H^0(\det (\EEE (1)))=
H^0 (\OOO (2))$ and that $\mt{ev}$ is an isomorphism.
This is proved in the next lemmas. The uniqueness modulo
$PGL(6)$ is proved in Lemma \ref{PGL6}.

\end{proof}

\begin{lemma}
Let $\mt{Pf}_2:\PPPP\dasharrow\PPPP$ be the Pfaffian map, sending a
skew-symmetric $6\times 6$ matrix $M$ to the collection
of its 15 quadratic Pfaffians. Then $\mt{Pf}_2^2=
\id_\PPPP$, the restriction of $\mt{Pf}_2$ to $\PPPP\setminus\Xi$
is an isomorphism onto $\PPPP\setminus G$, and $i=\mt{Pf}_2|_X$.
\end{lemma}

Thus $\mt{Pf}_2$ is an example of a Cremona quadratic transformation.
Such transformations were studied in \cite{E-SB}.

\begin{proof}
Let $(e_i)$, $(\eps_i)$ be dual bases of $E,E^\dual$ respectively,
and $(e_{ij}=e_i\wedge e_j),(\eps_{ij})$ the corresponding bases
of $\wedge^2E$, $\wedge^2E^\dual$. Identify $M$ in the source of
$\mt{Pf}_2$ with a 2-form $M=\sum a_{ij}\eps_{ij}$. Then $\mt{Pf}_2$
can be given by the formula
$\mt{Pf}_2(M)=\frac{1}{2!4!}M\wedge M\inner e_{123456}$,
where $e_{123456}=e_1\wedge\ldots\wedge e_6$, and $\inner$
stands for the contraction of tensors. Notice that $\mt{Pf}_2$
sends 2-forms of rank 6,4, resp. 2 to bivectors of rank 6,2, resp. 0.
Hence $\mt{Pf}_2$ is not defined on $G'$ and contracts $\Xi '
\setminus G'$ into $G$. In fact, the Pfaffians of a 2-form $M$
of rank 4 are exactly the Pl\"ucker coordinates of $\ker M$,
which implies $i=\mt{Pf}_2|_X$.

In order to iterate $\mt{Pf}_2$, we have to identify its source
$\PP (\wedge^2E^\dual )$ with its target $\PP (\wedge^2E )$.
We do it in using the above bases: $\eps_{ij}\mapsto e_{ij}$.
Let $N=\mt{Pf}_2^2(M)=\sum b_{ij}\eps_{ij}$.
Then each matrix element $b_{ij}=b_{ij}(M)$ is a polynomial of
degree 4 in $(a_{kl})$, vanishing on $\Xi '$. Hence it is divisible
by the equation of $\Xi '$, which is the cubic Pfaffian $\mt{Pf}(M)$.
We can write $b_{ij}=\tilde{b}_{ij}\mt{Pf}(M)$, where $\tilde{b}_{ij}$
are some linear forms in $(a_{kl})$. Testing them on a collection of
simple matrices with only one variable matrix element, we find
the answer: $\mt{Pf}_2(M)=\mt{Pf}(M)M$. Hence $\mt{Pf}_2$ is a birational
involution.
\end{proof}

\begin{lemma}\label{iofaline}
Let $l\subset V_3$ be a line. Then $i(l)$ is a conic in $\PPPP$,
and the lines of $\PP^5$ parametrized by the points of $i(l)$
sweep out a quadric surface of rank $3$ or $4$.
\end{lemma}

\begin{proof}
The   restriction of $\mt{Cl}$ to the lines
in $V_3$ is written out in \cite{AR} on pages 170 (for a
non-jumping line of $\KKK$, formula (49.5)) and 171 (for a jumping line).
These formulas imply the assertion; in fact, the quadric surface
has rank 4 for a non-jumping line, and rank 3 for a jumping one.
\end{proof}

\begin{lemma}\label{ARlem}
The map $i$ is injective.
\end{lemma}

\begin{proof}
Let  $\tilde{\Xi}$ be the natural desingularization of $\Xi '$ 
para\-metrizing
pairs $(M,l)$, where $M$ is a skew-symmetric $6\times 6$
matrix and $l$ is a line in the projectivized kernel
of $M$. We have $\tilde{\Xi}=\PP (\wedge^2(E_X/\SSS))$,
where $\SSS$ is the tautological rank 2 vector bundle on $G=G(2,6)$.
$\tilde{\Xi}$ has two natural projections 
$p:\tilde{\Xi}\lra G\subset \PPPP$ and
$q:\tilde{\Xi}\lra \Xi '\subset \PP^{14\dual}$. The classifying map
of $\KKK$ is just $\mt{Cl}=pq^{-1}$. $q$
is isomorphic over the alternating forms of rank 4, so
$q^{-1}(V_3)\simeq V_3$. $p$ is at least bijective on $q^{-1}(V_3)$.
In fact, it is easy to see that the fibers of $p$ can only be
linear subspaces of $\PPPP$. Indeed, the fiber of $p$ is nothing but
the family of matrices $M$ whose kernel contains a fixed plane, 
hence it is a linear subspace $\PP^5$ of $\PPPP^\dual$, and
the fibers of $p|_{q^{-1}(V_3)}$ are $\PP^5\cap V_3$.
As $V_3$ does not contain planes, the only possible fibers are
points or lines. By the previous lemma, they can be only points,
so $i$ is injective.
\end{proof}

\begin{lemma}
$i$ is defined by  the image of 
the map $\mt{ev}:\Lambda^2H^0(\EEE (1))$
$\lra $
$H^0(\det (\EEE (1)))=$
$H^0 (\OOO (2))$ considered as a linear subsystem of $|\OOO (2)|$.
\end{lemma}

\begin{proof}
Let $(x_i=\eps_i)$ be the coordinate functions on $E$, dual to the
basis $(e_i)$. The $x_i$ can be considered
as sections of $\KKK^\dual$. Then $x_i\wedge x_j$ can be considered
either as an element $x_{ij}$ of $\wedge^2E^\dual=\wedge^2H^0(\KKK^\dual)$, or
as a section $s_{ij}$ of $\wedge^2\KKK^\dual$. For a point $x\in V_3$, the
Pl\"ucker coordinates of the corresponding plane $K_x\subset E$
are $x_{ij}(\nu )$ for a non zero bivector $\nu\in\wedge^2K_x$.
By construction, this is the same as $s_{ij}(x)(\nu )$. This proves
the assertion.

\end{proof}

\begin{lemma}\label{PGL6}
Let $X\tilde{\lra}\Xi '\cap\La_1$,  $X\tilde{\lra}\Xi '\cap\La_2$
be two representations of $X$ as linear sections of $\Xi '$,
$\KKK_1,\KKK_2$ the corresponding kernel bundles on $X$.
Assume that $\KKK_1\simeq\KKK_2$. Then there exists a linear
transformation $A\in GL(E^\dual )=GL_6$ such that 
$\Xi '\cap\wedge^2A(\La_1)$ and $\Xi '\cap\La_2$ have the same image under
the classifying maps into $G$. The family of linear sections $\Xi '\cap\La$
of the Pfaffian cubic with the same image in $G$ is a rationally 1-connected
subvariety of $G(5,15)$, generically of dimension~$0$.
\end{lemma}

\begin{proof}
The representations  $X\tilde{\lra}\Xi '\cap\La_1$,  
$X\tilde{\lra}\Xi '\cap\La_2$ define two isomorphisms
$f_1:E^\dual\lra H^0(\KKK_1)$, $f_2:E^\dual\lra H^0(\KKK_2)$.
Identifying $\KKK_1,\KKK_2$, define $A=f_2^{-1}\circ f_1$.

Assume that $\La =\wedge^2A(\La_1)\neq\La_2$. Then the two
3-dimensional cubics $\Xi '\cap\La$ and $\Xi '\cap\La_2$ are isomorphic
by virtue of the map $f=f_2\circ f_1^{-1}\circ (\wedge^2A)^{-1}$.
By construction, we have $\ker M=\ker f(M)$ for any $M\in \Xi '\cap\La$.
Hence $\Xi '\cap\La$ and $\Xi '\cap\La_2$ represent two cross-sections of
the map $pq^{-1}$ defined in the proof of Lemma \ref{ARlem} over
their common image $Y=pq^{-1}(\Xi '\cap\La )=pq^{-1}(\Xi '\cap\La_2 )$,
and $f$ is a morphism over $Y$.
These cross-sections do not meet the
indeterminacy locus $G'\subset\Xi '$ of $pq^{-1}$, because it is
at the same time the singular locus of $\Xi '$ and both 3-dimensional cubics
are nonsingular.
The fibers of $pq^{-1}$ being linear subspaces of $\PPPP ^\dual$, the generic
element of a linear pencil $X_{\la :\mu}=\Xi '\cap (\la \La +\mu \La_2)$
represents also a cross-section of $pq^{-1}$ that does not meet  $G'$.
So there is a one-dimensional family of representations of $X$ as
a linear section of the Pfaffian cubic which are not equivalent under
the action of $PGL(6)$ but induce the same vector bundle $\KKK$.
This family joins $\Xi '\cap\La$ and $\Xi '\cap\La_2$ and its
base is an open subset of $\PP^1$.
This cannot happen for generic $X,\EEE$, because both the family
of vector bundles $\EEE$ and that of representations of $X$ as a linear section
of $\Xi '$ are 5 dimensional for generic $X$ (Theorem \ref{mainMT}
and Proposition \ref{puts}).

\end{proof}

\begin{lemma}
For a generic $3$-dimensional linear section $V_3$ of $\Xi '$,
the $15$ quadratic Pfaffians of $M\in V_3$ are linearly
independent in $|\OOO_{V_3}(2)|$.
\end{lemma}

The authors of \cite{IR} solve a similar problem: they
describe the structure of the restriction of  $\mt{Pf}_2$
to a 4-dimensional linear section of the Pfaffian cubic.

\begin{proof}
It is sufficient to verify this property for a special $V_3$.
Take Klein's cubic
$$
v^2w+w^2x+x^2y+y^2z+z^2v=0.
$$
Adler (\cite{AR}, Lemma (47.2)) gives the representation
of this cubic as the Pfaffian of the following matrix:
$$ M=\left[\begin{array}{cccccc} 0&v&w&x&y&z\\
-v&0&0&z&-x&0\\
-w&0&0&0&v&-y\\
-x&-z&0&0&0&w\\
-y&x&-v&0&0&0\\
-z&0&y&-w&0&0\\ \end{array}\right]
.$$

Its quadratic Pfaffians are given by
$$
c_{ij}=(-1)^{i+j+1}(a_{pq}a_{rs}-a_{pr}a_{qs}+a_{ps}a_{qr}),
$$
where $p<q<r<s$, $(pqrsij)$ is a permutation of $(123456)$, and $(-1)^{i+j+1}$
is nothing but its sign.
A direct computation shows that the 15 quadratic Pfaffians are linearly independent.
\end{proof}

This ends the proof of Theorem \ref{KKK}.

\section{Minimal sections of 2-dimensional conic bundle}

Let $X$ be a generic cubic threefold.
To prove the irreducibility of the fibers of the Abel-Jacobi map $\Phi$
of Theorem \ref{mainMT}, we will  use other Abel--Jacobi maps.
Let us fix a line $l_0$ in $X$, and denote by $\Phi_{d,g}$  the Abel--Jacobi map
of the family $H_{d,g}$ of curves of degree $d$ and of
arithmetic genus $g$ in $X$
having $dl_0$ as reference curve. The precise domain of definition
of $\Phi_{d,g}$ will be specified in the context in each particular
case. So, $\Phi_{5,1}$  will be exactly the above map $\Phi$ defined
on $\HHH$.

We will provide a description of $\Phi_{4,0}$, obtained by an 
application of the techiniques of \cite{Ili}. This map is
defined on the family of normal rational quartics in $X$. For completeness,
we will mention a similar description of $\Phi_{3,0}$, the Abel--Jacobi
map of twisted rational cubics in $X$. As was proved in \cite{MT}, these
families of curves are irreducible for a generic $X$. 

Let $L_0\subset X$ be a generic line, $p :\tilde{X}\lra\PP^2$
the projection from $L_0$, giving to $\tilde{X}=\Blowup_{L_0}(X)$
a structure of a conic bundle. Let $C\subset\PP^2$ be a generic conic, then
$Y=p^{-1}(C)$ is a 2-dimensional conic bundle, and $p_Y=p |_{Y}:Y\lra C$ is
the conic bundle structure map. It is well known (see \cite{B}),
that the discriminant
curve $\Delta\subset\PP^2$ of $p$ is a smooth quintic, and the components
of the reducible conics $\PP^1\vee\PP^1$ over points of $\Delta$ are
parametrized by a non-ramified two-sheeted covering $\pi : \tilde{\Delta}
\lra\Delta$. As $C$ is generic, there are 10 distinct points in $\Delta
\cap C$, giving us 10 pairs of lines $\{ l_1\cup l'_1\cup\ldots \cup l_{10}\cup
l'_{10}\} =p^{-1}(\Delta\cap C)$.
We will identify the components $l$
of reducible fibers of $p$ with points of $\tilde{\Delta}$,
so that $\{ l_1, l'_1,\ldots ,l_{10},
l'_{10}\}=\pi^{-1}(\Delta\cap C)\subset\tilde{\Delta}$.
Let $p_\alp :Y_\alp \lra C$ be any of the $2^{10}$ ruled surfaces 
obtained by contracting
the $l'_i$ with $i\in\alp$ and the $l_j$ with $j\not\in\alp$, where $\alp$
runs over the subsets of $\{1,2,\ldots ,10\}$. Then the $Y_\alp$
are divided into two classes: even and odd surfaces, according to the
parity of the integer $n\geq 0$ such that $Y_\alp\simeq\FF_n=\PP 
(\OOO_{\PP^1}\oplus\OOO_{\PP^1}(-n))$. Remark, that the surfaces
$Y_\alp$ are in a natural one-to-one correspondence with effective divisors $D$
of degree 10 on $\tilde{\Delta}$ such that $\pi_*D=\Delta\cap C$.
The 10 points of such a divisor correspond to lines ($l_i$ or $l'_i$)
which are not contracted by the map $Y\lra Y_\alp$. 
For a surface $Y_\alp$,
associated to an effective divisor $D$ of degree 10,
we will use the alternative notation $Y_D$.

The next theorem is a particular case of the result of \cite{Ili}. 

\begin{theorem}\label{minsec}
Let $X$ be a generic cubic threefold, $C\in\PP^2$ a generic conic.
Then, in the above notation, the following assertions hold:

(i) There are only two isomorphism classes of surfaces among the $Y_\alp$:\
$Y_{odd}\simeq\FF_1$ and $Y_{even}\simeq\FF_0\simeq\PP^1\times\PP^1$.

(ii) The family $\CCC_-$ of the proper transforms in $X$
of $(-1)$-curves in each one of the odd surfaces $Y_{\alp}\simeq \FF_1$ over
all sufficiently generic conics $C\subset\PP^2$ is identified
with a dense open subset in the family of twisted rational cubic
curves $C^3\in X$ meeting $L_0$ at one point. 

(iii) Let $\Phi_{3,0}$ be the Abel--Jacobi map of the family of rational
twisted cubics. Let $\Phi_-
=\Phi_3|_{\CCC_-}$ be its restriction.  Then $\Phi_-$ is onto an open 
subset of the theta divisor
of $J^2(X)$. For generic $C^3\in\CCC_-$, which is a proper transform
of the $(-1)$-curve in the ruled surface $Y_\alp$ associated to
an effective divisor $D_\alp$ of degree 10 on $\tilde{\Delta}$,
the fiber $\Phi_-^{-1}\Phi_-(C^3)$ can be identified with 
an open subset of $\PP^1=|D_\alp |$
by the following rule:
$$
D\in |D_\alp |\ \mapsto \ \left|\begin{minipage}{7 truecm}
{\scriptsize the proper transform in $X$
of the $(-1)$-curve in $Y_D$ if $Y_D\simeq\FF_1$}
\end{minipage}\right.
$$

(iv) Let $\CCC_+$ be the family of the  proper 
transforms in $X$ of the curves in the second ruling on any one of the
even surfaces $Y_\alp\simeq\PP^1\times\PP^1$
for all sufficiently generic conics $C$;
the second ruling means the one which is different from that
consisting of fibers of $\pi_\alp$. Then $\CCC_+$ is identified
with a dense open subset in the family of normal rational quartic
curves $C^4\in X$ meeting $L_0$ at two points.

(v) Let $\Phi_{4,0}$ be the Abel--Jacobi map of the family of
rational normal quartics. Let $\Phi_+
=\Phi_{4,0}|_{\CCC_+}$ be its restriction.  Then $\Phi_+$ is onto an open 
subset
of $J^2(X)$. For generic $C^4\in\CCC_+$ which is the proper transform
of a curve on the ruled surface $Y_\alp$ associated to
an effective divisor $D_\alp$ of degree 10 on $\tilde{\Delta}$,
where ${\dim |D_\alp |=0}$ and
the fiber $\Phi_+^{-1}\Phi_+(C^4)\simeq\PP^1$ consists of
the proper transforms of all the curves of the second ruling
on $Y_\alp$.
\end{theorem}

The irreducibility of $\Phi_+^{-1}\Phi_+(C^4)$ in the above statement
is an essential ingredient of the proof of the following theorem,
which is the main result of the paper.

\begin{theorem} \label{degree1}
Let $X$ be a nonsingular cubic threefold. Then
the degree of the \'etale map $\Psi$ from Theorem \ref{mainMT}
is $1$. Equivalently, all the fibers of the Abel--Jacobi map $\Phi$
are isomorphic to $\PP^5$.
\end{theorem}

This obviously implies:

\begin{corollary}\label{openinJ}
The open set $M\subset M_X(2;0,2)$ in the moduli space of
vector bundles on $X$ parametrizing those stable
rank 2 vector bundles which arise via Serre's construction
from normal elliptic quintics is isomorphic to an open
subset in the intermediate Jacobian of $X$.
\end{corollary}

We will start by the following lemma.

\begin{lemma} \label{Hliz}
Let $X$ be a generic cubic threefold.
Let $z \in J^2(X)$ be a generic point, $\HHH_i(z)\simeq{\bf P}^5$ 
any component of ${\Phi}^{-1}(z)$. Then, for any line $l \subset X_3$,
the family

$\HHH_{l;i}(z)$ :=
$\{ C \in \HHH_i(z) : C = l + C'$, where
$C'$ is a curve of degree~$4$~$\}$

\noindent is isomorphic to ${\bf P}^1$.
\end{lemma}

\begin{proof}
By Theorem \ref{mainMT}, the curve
$C$ represented by the generic point of $\HHH_{i}(z)$ is
a (smooth) normal elliptic quintic. Let $\EEE =\EEE_C$ be the
associated vector bundle, represented by the point $\phi ([C])\in M$.
Choose any representation of $X$ as a linear section of the
Pfaffian cubic $\Xi '$ as in Theorem \ref{KKK}, so that
$\EEE (1)\simeq\KKK^\dual$. The projective space $\HHH_i(z)$
is naturally identified with $\PP^{5\dual}=\PP (E^\dual )$. This follows from
the proof of Theorem \ref{KKK}. Indeed,
the curves $C$ represented by points of $\HHH_i(z)$ are exactly the
zero loci of the sections of $\EEE (1)$, and the latter are induced
by linear forms on $E$ via the natural surjection $E_X\lra\KKK^\dual$.
The zero loci of these sections are of the form $\Cl^{-1}(\si_{11}(L))$,
where $L\in\PP^{5\dual}$ runs over all the hyperplanes in $\PP^5$.

Let $l$ be a line in $X$. By Lemma \ref{iofaline}, the quadratic
pencil of lines with base $\mt{Cl}(l)$ sweeps out a quadric surface
$Q(l)$ of rank 3 or 4. Let $<\! Q(l)\! >\simeq\PP^3$ be the linear span of
$Q(l)$ in $\PP^5$. Then $l$ is a component of $\mt{Cl}^{-1}(\si_{11}(L))$
if and only if $<\! Q(l)\! >\subset L$. Such hyperplanes $L$ form the
pencil $<\! Q(l)\! >^\dual\simeq\PP^1$ in $\PP^{5\dual}$. Obviously,
the pencil $\{ \mt{Cl}^{-1}(\si_{11}(L))\mid L\in<\! Q(l)\! >^\dual\}$ contains
exactly all the curves, represented by points of $\HHH_i(z)$
and having $l$ as an irreducible component.
\end{proof}

Now our aim is to  show that the generic member of $\HHH_{l;i}(z)$ 
is a rational normal quartic having $l$ as one of its chords. 
Then we will be able to 
apply the description of such curves given by Theorem \ref{minsec}, (iv),~(v).

\section{Smoothing \( C'+l\)}\label{EEECl}

Let $X$ be a nonsingular cubic threefold, $C=C'+l\subset X$ 
a rational normal quartic plus one of its chords. Then one can apply
Serre's construction (\ref{serre}) to $C$
to obtain a self-dual rank 2 vector bundle $\EEE =\EEE_C$ in $M_X(2;0,2)$
like it was done in \cite{MT} for a nonsingular
$C$. One proves directly that $\EEE$ possesses all the essential
properties of the vector bundles constructed from normal
elliptic quintics. First of all, our $C$ is a locally complete
intersection in $X$ with trivial canonical sheaf $\omega_C$, and this
implies (see the proofs of Lemma 2.1 and Corollary 2.2 in loc. cit.) that
$\Ext^1(\III_C(2),\OOO_X )\simeq H^0(C,\omega_C)\simeq\CC$
and that $\EXT_{\OOO_X}^1(\III_C(2),\OOO_X)=
\EXT_{\OOO_X}^2(\OOO_C,\omega_X)=\omega_C$, so that $\EEE$
is uniquely determined up to isomorphism
and is locally free. One can also easily show that
$h^0(\III_C(1))=h^1(\III_C(1))=h^2(\III_C(1))=0$, and this implies
(see the proofs of Corollary 2.4, Proposition 2.6 and Lemma 2.8
in loc. cit.) the stability of $\EEE$
and the fact that the zero loci of nonproportional sections of
$\EEE (1)$ are distinct complete intersection linearly normal
quintic curves. Further, remark that
$h^0(\III_C (2))=5$  (the basis of $H^0(\III_C (2))$ is given
in appropriate coordinates in (\ref{Qi}) below);  the restriction exact sequence
\begin{equation}\label{restriction}
0\lra\III_C(k)\lra\OOO_X(k)\lra\OOO_C(k)\lra 0 
\end{equation}
with $k=2$ implies also $h^i(\III_C (2))=0$ for $i>0$.
One deduces from here
$h^0(\EEE (1))=6,h^i(\EEE (1))=0$ for $i>0$. Hence the sections
of $\EEE$ define a $\PP^5$ in $\Hilb$. 

We want to show that,
generically, this $\PP^5$ is of the form $\HHH_i(z)$, that is $\EEE$ has a
section whose zero locus is a (smooth) normal elliptic quintic.
First, we make a routine verification that $C$ can be smoothed
into a normal elliptic quintic, that is $[C]\in\Hilb$ is in the
closure of $\HHH^*$. Afterwards, we will show that the smoothing
can be effectuated inside the $\PP^5$ of zero loci of sections
of $\EEE$. To this end, we will find an {\em example}
of $(X,\EEE )$, in which the locus of the 
curves of type `normal rational quartic
plus its chord' inside the $\PP^5$ has at least one 3-dimensional component. 
By a standard dimension count,
this will imply that all the components of this locus are 
3-dimensional for generic $(X,\EEE )$, and that the generic point of
the $\PP^5$ represents a nonsingular curve.

\begin{lemma}\label{quartchord}
Let $X$ be a nonsingular cubic threefold, $C=C'+l\subset X$ 
a rational normal quartic plus one of its chords, $\EEE$ the
vector bundle defined in \ref{EEECl}. Then the following
assertions are true:

(i)  $h^i(\EEE (-1))=0\;\; \forall\; i\in\ZZ$, hence $\EEE$
is an instanton vector bundle of charge 2. Further,
$h^0(\EEE\otimes\EEE )=1$, $h^1(\EEE\otimes\EEE )=5,
h^2(\EEE\otimes\EEE )=h^3(\EEE\otimes\EEE )=0$, hence
$M_X(2;0,2)$ is smooth of dimension~$5$ at $[\EEE ]$.

ii) 
$h^0(\NNN_{C/X})=10,h^1(\NNN_{C/X})=0$, hence $\Hilb$ is
smooth of dimension 10 in $[C]$. Moreover, if we assume that $
\NNN_{C'/X}\not\simeq\OOO\oplus\OOO (6)$ or $\NNN_{l/X}\not\simeq
\OOO (-1)\oplus\OOO (1)$, then $C$ is strongly smoothable and a sufficiently
small deformation ${\mathfrak C}\lra U$
of $C$ parametrizes curves of only the following three types:
(a) for $u$ in a dense open subset of $U$, $C_u$ is a normal elliptic
quintic; (b) over on open subset of a divisor $\Delta_1\subset U$, $C_u$ is a linearly normal
rational curve with only one node as singularity;
(c)~over a closed subvariety of pure codimension 2 $\Delta_2\subset U$, $C_u$ is 
of the same type as $C$, that is a
normal rational quartic plus one of its chords.
\end{lemma}

\begin{proof}
As concerns the numerical values for the $h^i$, the proof goes exactly
as that of Lemma 2.7 in \cite{MT} with only one modification: the authors used
there the property of a normal elliptic quintic $h^0(\NNN_{C/\PP^4}(-2))$
$=$ $0$,
proved in Proposition V.2.1 of \cite{Hu}. Here we should verify directly this
property for our curve $C=C'+l$. This is an easy exercise: one can use
the identifications of the normal bundles of $C',l$
\begin{equation}\label{normquartic}
\NNN_{C'/{\PP^4}}\simeq 3\OOO_{\PP^1}(6)\ ,\qquad \NNN_{l/{\PP^4}}
\simeq 3\OOO_{\PP^1}(1)
\end{equation}
and the three natural exact sequences
\begin{equation}\label{seq1}
0\lra\NNN_{C/{W}}\lra \NNN_{C/{W}}|_{C'}\oplus \NNN_{C/{W}}|_l
\lra \NNN_{C/{W}}\otimes\CC_S\lra 0,
\end{equation}
\begin{equation}\label{seq2}
0\lra \NNN_{C'/{W}}\lra \NNN_{C/{W}}|_{C'}\lra T^1_S\lra 0,
\end{equation}
\begin{equation}\label{seq3}
0\lra \NNN_{l/{W}}\lra \NNN_{C/{W}}|_l\lra T^1_S\lra 0,
\end{equation}
where $S=\{ P_1,P_2\} =C'\cap l$, $\CC_S=\CC_{P_1}\oplus\CC_{P_2}$
is the sky-scraper sheaf
with the only nonzero stalks at $P_1,P_2$ equal to $\CC$,  $W=\PP^4$,
and $T^1_{S}$ denotes
Schlesinger's $T^1$ of a singularity; we have $T^1_{S}\simeq\CC_S$
for nodal curves.

The values of $h^i(\EEE\otimes\EEE )$ in (i) imply the stated
properties of the moduli space, because $\EEE$ is self-dual, and
so $h^i(\EEE\otimes\EEE )=\dim\Ext^i(\EEE ,\EEE )$.

For the remaining assertions of (ii), we will apply Theorem 4.1 of \cite{HH}
\footnote{Hartshorne--Hirschowitz formulated it for nodal curves in $\PP^3$, but the
proof and the techniques of the paper remain valid if one replaces
$\PP^3$ by any nonsingular projective variety; see Remark 4.1.1 in \cite{HH}.}.
It states that if the elementary transformations of the normal
bundles to $C', l$ satisfy
$H^1(C',\mt{elm}_{P_i}^+\NNN_{l/X})=
H^1(l,\mt{elm}_{S}^-\NNN_{C'/X})=0$, then
$C$ is strongly smoothable. In fact, by (\ref{normquartic})
and
$$
0\lra \NNN_{C'/X}\lra\NNN_{C/{\PP^4}}\lra \NNN_{X/{\PP^4}}|_C\lra 0\ ,
$$
we see that $\NNN_{C'/X}\simeq\OOO (a)\oplus\OOO (b)$ with
$a+b=6, 0\leq a\leq b\leq 6$. By \cite{CG}, $\NNN_{l/X}\simeq 2\OOO$
or $\OOO (-1)\oplus\OOO (1)$, so $\mt{elm}_{P_i}^+\NNN_{l/X}$
$=$ $\OOO (-1)\oplus\OOO (2)$ or
$\OOO \oplus\OOO (1)$.
For $C'$, $\mt{elm}_{S}^-\NNN_{C'/X}$ may be one of the sheaves
$\OOO (a-2)\oplus\OOO (b)$, $\OOO (a-1)\oplus\OOO (b-1)$,
or $\OOO (a)\oplus\OOO (b-2)$.  So, the hypotheses
of the theorem may be not verified only if $a=0,b=6$. In interchanging
the roles of $C',l$ and assuming that $\NNN_{l/X}\simeq 2\OOO$, we can see that 
$\mt{elm}_{S}^-\NNN_{l/X}\simeq 2\OOO (-1)$ or $\OOO\oplus\OOO (-2)$. The second
case is impossible, because $l$ and the tangent directions of $C'$ at the points
$P_i$ are not coplanar, so the centers of the elementary transformation
$\tilde{P}_i\in\PP (\NNN_{l/X}|_{P_i})\simeq\PP^1$, corresponding to
the directions of $C'$ at $P_i$, do not lie on the
same section of $\PP (\NNN_{l/X})=\PP^1\times\PP^1$.

Thus, in both cases the theorem can be applied, and we conclude that the natural maps
$\delta_i: H^0(\NNN_{C/X})\lra T^1_{P_i}C=\CC_{P_i}$ are surjective. Hence the discriminant
divisor $\Delta_1\subset U$ has locally analytically two nonsingular
branches with tangent spaces $\ker\delta_i\subset H^0(\NNN_{C/X})=T_{[C]}\Hilb$,
each unfolding only one of the two singular points of $C$,
and their transversal intersection $\Delta_2$ parametrizes the deformations
preserving the two singular points. To conclude the proof, remark that
the linear normality and $p_a(C)$ are preserved under small deformations.
\end{proof}

\begin{lemma}\label{smoothHE}
Let $X$ be a generic cubic threefold, $C=C'+l\subset X$ 
a generic rational normal quartic plus one of its chords, $\EEE$ the
vector bundle defined above. Let $\HHH_\EEE\subset\Hilb$
be the $\PP^5$ of zero loci of sections of $\EEE (1)$.
Then the assumption of Lemma \ref{quartchord}, (ii) 
for the normal bundle of $l$ is verified and, moreover,
$\dim\Delta_i\cap\HHH_\EEE =5-i$ for $i=1,2$. This implies that
$C'+l$ can be smoothed not only inside $\Hilb$, but also inside  $\HHH_\EEE$.
\end{lemma}

\begin{proof}
We have to show that $\NNN_{l/X}\simeq 2\OOO$ and that the
natural map $T_{[C]}\Hilb\lra T^1_SC$ \hfill remains surjective when
restricted to \hfill $T_{[C]}\HHH_\EEE$ \hfill $\subset $ \linebreak
$
T_{[C]}\Hilb$. It suffices
to do this only for one special cubic threefold $X$ and for
one special $C$, because both conditions are open. 
So, choose a curve $C$ of type $C'+l$ in $\PP^4$,
then a cubic $X$ passing through $C$. Take, for example, the closures
of the following affine curves:
$$
C'\  =\{ x_1=t,\ x_2=t^2,\ x_3=t^3,\ x_4=t^4 \}\ ,\ \
l\ =\{ x_1=x_2=x_3=0 \}\ .
$$
The family of quadrics passing through $C$ is 5-dimensional with generators
\begin{multline}\label{Qi}
Q_1=x_2-x_1^2,\ Q_2=x_3-x_1x_2,\ Q_3=x_1x_3-x_2^2,\\
Q_4=x_1x_4-x_2x_3,\ Q_5=x_2x_4-x_3^2\ .
\end{multline}
The cubic hypersurface in $\PP^4$ with equation $\sum\alpha_i(x)Q_i$ is
nonsingular for generic linear forms $\alpha_i(x)$, so we can choose $X$
to be of this form. We verified, in using the Macaulay program \cite{BS},
that the choice $\alp_1= 0,\alp_2=-1,\alp_3=x_2,
\alp_4=-x_1,\alp_5=x_4$ yields a nonsingular
$X=\{ x_1x_2-x_2^3-x_3+2x_1x_2x_3-x_1^2x_4-x_3^2x_4+x_2x_4^2=0\}$ such that $\NNN_{l/X}\simeq 2\OOO$.

Look at the following commutative diagram with exact rows
and columns, where the first row is the restriction of (\ref{serre})
to the subsheaf of the sections of $\EEE (1)$ vanishing along $C$.

\begin{equation}\label{CD2}
\begin{CD}
 @. @. 0 @. 0 @. \\
@. @. @VVV @VVV @. \\
0 @>>> \OOO_X @>>> \EEE (1)\otimes\III_C @>>> \III_C^2(2) @>>> 0 \\
@. @| @VVV @VVV @. \\
0 @>>> \OOO_X @>>> \EEE (1) @>>> \III_C(2) @>>>  0 \\
@. @. @VVV @VVV @. \\
 @.  @. \NNN_{C/X} @= \NNN^\dual _{C/X}(2) @.  \\
@. @. @VVV @VVV @. \\
 @. @. 0 @. 0 @. \\
\end{CD}
\end{equation}

It allows to identify the tangent space $T_{[\EEE ]}\HHH_\EEE=
H^0(\EEE (1))/H^0(\EEE (1)\otimes\III_C)$ with the image of $H^0(\III_C(2))$
in $H^0(\NNN_{C/X})=H^0(\III_C(2)/\III_C(2)^2)$. So, we have to show that
the derivative $d:H^0( \III_C(2))\lra T^1_SC$ is surjective. Using the
basis (\ref{Qi}) of $H^0(\III_C(2))$, we easily verify that this is
the case (in fact, $dQ_1,dQ_2$ generate $T^1_SC$).
\end{proof}

\begin{lemma}\label{dim3}
With the hypotheses of Lemma \ref{Hliz}, 
the family $\CCC_i(z)$ of curves of the form $C'+l$ in $\HHH_i(z)$, where
$C'$ is a rational normal quartic and $l$ one of its chords,
is non-empty and equidimensional of dimension~$3$.
\end{lemma}

\begin{proof}
According to \cite{MT}, the family of rational normal quartics
in a nonsingular cubic threefold $X$ has dimension 8, and is irreducible
for generic $X$. By Theorem \ref{TT},
each rational normal quartic $C'$ has exactly 16 chords $l$ in $X$,
so the family $\Delta_2 =\Delta_2 (X)$ of pairs $C'+l$ is equidimensional 
of dimension 8.
It suffices to verify that one of the components of $\Delta_2$, say 
$\Delta_{2,0}$,
meets $\HHH_i(z)$ at some point $b$ with local
dimension $\dim_b\Delta_{2,0}\cap\HHH_i(z)=3$
for one special cubic threefold $X$, for
one special $z$ and for at least one $i$. But this was done
in the previous lemma. Indeed, the fact that $C$ can be smoothed inside
$\HHH_\EEE$ implies that $\EEE\in\HHH$, hence $\HHH_\EEE=\HHH_i(z)$
for some $i,z$.
The assertion for general $X,z$
follows by the relativization over the family of cubic threefolds
and the standard count of dimensions. 
\end{proof}

\begin{corollary}\label{rnqpc}
With the hypotheses of Lemma \ref{Hliz}, let $l$ be a generic
line in $X$. Then the generic member of the pencil $\HHH_{l;i}(z)$
is a rational normal quartic plus one of its chords.
\end{corollary}

\begin{proof}
We know already that the family of pairs $C'+l\in \CCC_i(z)$
is 3-dimensional. Now we are to show that the second components
$l$ of these pairs move in a dense open subset in the Fano surface $F$
of $X$. This can be done by an infinitesimal argument: let $C=C'+l$
be such a pair. Then a small neighborhood of $[C]$ in $\HHH_i(z)$
contains a smooth subvariety $D$ of codimension 2 with tangent
space $T_{[C]}D=\ker\{ T_{[C]}\HHH_i(z)\lra T^1_SC\}$, which
parametrizes the curves $C''+l$ of the same type. It suffices to show
that the natural projection of $T_{[C]}D$ to $T_{[l]}F=H^0(\NNN_{l/X})$ is surjective.

The exact triples (\ref{seq1}), (\ref{seq2}), (\ref{seq3}) with $W=X$
together with the observation that all the $H^1$'s vanish imply that
the natural map $H^0(\NNN_{C/X})$ $\lra $ $H^0(\NNN_{C/X}|_l)$
is surjective and that it restricts to a surjective map between
the kernels of the respective maps to $T^1_SC$:
\begin{multline*}
T_{[C]}\Delta_2X=\ker\{ H^0(\NNN_{C/X})\lra T^1_SC\}\surj\\
H^0(\NNN_{l/X})=\ker\{ H^0(\NNN_{C/X}|_l)\lra T^1_SC\} .
\end{multline*}
We want to see that it will remain surjective even if we shrink
its source to $T_CD\subset T_{[C]}\Delta_2X$.

Let $\EEE =\EEE_C$ be as above. Then the triple (\ref{serre})
determines an isomorphism $\NNN_{C/X}=\EEE (1)|_C$ in such a way
that $T_{[C]}\HHH_i(z)$ is the image of the restriction map
$H^0(\EEE (1))\lra H^0(\EEE (1)|_C)$. Hence we can identify
$\NNN_{C/X}|_l$ with $\EEE (1)|_l$. Look at the exact triple
$$
0\lra\EEE (1)\otimes\III_l\lra\EEE (1)\lra\EEE (1)|_l\lra 0 .
$$
By Lemma \ref{Hliz}, $h^0(\EEE (1)\otimes\III_l)=2$. By Lemma
\ref{smoothHE}, we have $\NNN_{l/X}\simeq 2\OOO$ for generic $C'+l$.
As in the proof of Lemma \ref{quartchord}, the fact that the centers
of the elementary transformation $\mt{elm}^+_S\NNN_{l/X}\simeq
\NNN_{C/X}|_l$ of $\NNN_{l/X}$ do not lie on the
same section of $\PP (\NNN_{l/X})=\PP^1\times\PP^1$ implies that
$\NNN_{C/X}|_l\simeq 2\OOO (1)$ and that the map
$\NNN_{C/X}|_l\lra T^1_SC$ is surjective. As $h^0 (\EEE (1))=6$, the last
exact triple gives the surjectivity of $H^0(\EEE (1))\lra H^0(\NNN_{C/X}|_l)$.
Restricting the map to the kernels of the natural surjections onto
$T^1_SC$, we obtain the result.

\end{proof}

\section{Fibers of \( \Phi_{4,0},\Phi_{5,1}\) and periods of varieties
\( V_{14}\) }

Now we are able to prove Theorem \ref{degree1}.
Let $X$ be a generic cubic threefold. Let $\Phi_{1,0}$, 
$\Phi_{4,0}$, resp. $\Phi^* =\Phi^*_{5,1}$ be the Abel--Jacobi map of lines,
rational normal quartics, resp. elliptic
normal quintics.
We will use the notation $\Phi$, or
$\Phi_{5,1}$ for the extension of $\Phi^*$ defined in the statement of
Theorem \ref{mainMT}. By Lemma \ref{smoothHE}, the generic curves of the
form $C'+l$, where $C'$ is a rational normal quartic and $l$ one
of its chords, are elements of $\HHH$, the domain of $\Phi$.

\smallskip

{\em Proof of Theorem \ref{degree1}.}
Let $z \in J^2(X)$ be a generic point, $\HHH_i(z)\simeq{\bf P}^5$ 
any component of ${\Phi}^{-1}(z)$. Choose a generic line $l$ on $X$.
In the notations of Lemma \ref{Hliz},
the number of pencils $\HHH_{l;i}(z)\simeq \PP^1$ with
generic member $C'_i+l$, where $C'_i$ is a rational normal quartic meeting
$l$ quasi-transversely at 2 points, and mapped to the same point $z$
of the intermediate Jacobian, is equal to the degree $d$ of $\Psi$.
Now look at the images of the curves $C'_i$ arising in these pencils
under the Abel--Jacobi map $\Phi_{4,0}$. Denoting $AJ$ the Abel--Jacobi
map on the algebraic 1-cycles homologous to 0, we have $AJ((C'_i+l)-
(C'_j+l))=AJ(C'_i-C'_j)=z-z=0$. Hence $\Phi_{4,0}(C'_i)=\Phi_{4,0}(C'_j)$
is a constant point $z'\in J^2(X)$. According to Theorem \ref{minsec},
the family of the normal rational quartics in a generic fiber of $\Phi_{4,0}$
meeting a generic line at two points is irreducible and is parametrized by
(an open subset of) a $\PP^1$. The point $z'$ is a generic one,
because $\Phi_{4,0}$ is dominant, and every rational normal quartic
has at least one chord. Hence $d=1$ and we are done.
\hfill\square

\begin{corollary}
$M,\HHH$ are irreducible and the degree of $\Psi$ is 1 not only
for a generic cubic $X$, but also for every nonsingular one.
\end{corollary}

\begin{proof}
One can easily relativize the constructions of $\HHH ,M,\Phi, \phi , \Psi$, etc.
over a small analytic (or \'etale) connected open set $U$ in the parameter space
$\PP^{34}$ of 3-dimensional cubics, over which all the cubics $X_u$ are nonsingular. 
We have to restrict ourselves to a
``small" open set, because we need a local section of the family
$\{\HHH_u\}$ in order to define the maps $\Phi, \Psi$.

The fibers $\HHH_u, M_u$  are equidimensional
and nonsingular of dimensions $10, 5$ respectively. Moreover, it is easy
to see that a normal elliptic quintic $C_0$ in a special fiber $X_{u_0}$
can be deformed to the neighboring fibers $X_u$. Indeed, one can embed
the pencil $\lambda X_{u_0}+\mu X_u$ into the linear system of hyperplane
sections of a 4-dimensional cubic $Y$ and show that the local dimension
of the Hilbert scheme of $Y$ at $[C_0]$ is 15, which implies that $C_0$
deforms to all the nearby and hence to all the nonsingular hypeplane sections
of $Y$.

Hence the families $\{\HHH_u\} ,\{ M_u\}$ are irreducible, flat of relative
dimensions 10, resp. 5 over $U$, and the degree of $\Psi$ is constant over $U$.
If there is a reducible fiber $M_u$, then the degree sums up over
its irreducible components, so it has to be strictly greater than 1.
But we know, that $d$ is 1 over the generic fiber, hence all the fibers
are irreducible and $d=1$ for all $u$.
\end{proof}

We are going to relate the Abel--Jacobi mapping of elliptic normal
quintics with that of rational normal quartics.
With our convention for the choice of reference curves in the
form $dl_0$ for a line $l_0$, fixed once and forever, we
have the identity
$$
\Phi_{5,1} (C'+l)= \Phi_{4,0}(C')+\Phi_{1,0}(l)\ .
$$

\begin{theorem}\label{fiberquart}
Let $X$ be a generic cubic threefold, $z \in J^2(X)$ a generic point.
Then the corresponding
fiber ${\Phi}_{4,0}^{-1}(z)$ is an irreducible
nonsingular variety of dimension~$3$, birationally equivalent
to $X$.
\end{theorem}

\begin{proof}
As we have already mentioned in the proof of Lemma \ref{dim3},
the family $H_{4,0}$ of rational normal quartics in $X$ is
irreducible of dimension 8. The nonsingularity of $H_{4,0}$
follows from the evaluation of the normal bundle of a rational
normal quartic in the proof of Lemma \ref{quartchord}. 
We saw also that $\Phi_{4,0}:H_{4,0}\lra J^2(X)$ is
dominant, so the generic fiber is equidimensional of dimension 3 and we have
to prove its irreducibility.

Let $\pi :\tilde{U}\lra U$ be the quasi-finite covering of $U=
\Phi (\HHH )$ parametrizing the irreducible components of the
fibers of $\Phi_{4,0}$ over points of $U$. Let $z\in U$ be generic,
and $\HHH_z\simeq\PP^5$ the fiber of $\Phi$. By Corollary \ref{rnqpc},
for a generic line $l$, 
we can represent $z$ as $\Phi_{4,0}(C')+\Phi_{1,0}(l)$ for
a rational normal quartic $C'$ having $l$ as one of its chords.
Let $\kappa :U\dasharrow \tilde{U}$ be the rational map sending $z$
to the component of $\Phi_{4,0}^{-1}\Phi_{4,0}(C')$ containing
$C'$. Let $\lambda =\pi\circ\kappa$. Theorem \ref{minsec} implies that $\lambda$
is dominant. Hence it is generically finite. Then
$\kappa$ is also generically finite, and we have for their
degrees $\deg\lambda =(\deg\pi )(\deg \kappa )$.

Let us show that $\deg\lambda =1$.
Let $z,z'$ be two distinct points in a generic fiber of $\lambda$.
By Theorem \ref{minsec}, $\Phi_{4,0}^{-1}\Phi_{4,0}(C')$
contains only one pencil of curves of type $C''+l$, where
$l$ is a fixed chord of $C'$, and $C''$ is a rational normal
quartic meeting $l$ in 2 points. But Lemma 
\ref{Hliz}  and Corollary \ref{rnqpc} imply that both
$\HHH_z$ and $\HHH_{z'}$ contain such a pencil. This
is a contradiction. Hence $\deg\lambda =\deg\pi =\deg\kappa =1$.

Now, choose a generic rational normal quartic $C'$
in $X$. We are going to show that $\Phi_{4,0}^{-1}\Phi_{4,0}(C')$
is birational to some $V_{14}$, associated to $X$, and hence birational to $X$
itself. Namely, take the  $V_{14}$ obtained by the Tregub--Takeuchi
transformation $\chi$ from $X$ with center $C'$. Let $x\in V_{14}$ be the
indeterminacy point of $\chi^{-1}$. The pair $(x,V_{14})$ is determined
by $(C',X)$ uniquely up to isomorphism, because $V_{14}$ is the image of
$X$ under the map defined by the linear system $|\OOO_X (8)-5C' |$
and $x$ is the image of the unique divisor of the linear system
$|\OOO_X(3-2C')|$.

By Theorem \ref{TT}, (iii), a generic $\xi\in V_{14}$ defines an inverse map
of Tregub--Takeuchi type from $V_{14}$ to the same cubic $X$.
As $X$ is generic, it has no biregular automorphisms, and hence this map defines
a rational normal quartic $\G$ in $X$. We obtain the rational map $\alp :
V_{14}\dasharrow H_{4,0}$, $\xi\mapsto [\G ]$, whose image contains $[C']$.
As $h^{1,0}(V_{14})=0$, the whole image $\alp (V_{14})$ is contracted to a point
by the Abel--Jacobi map. Hence, to show that $\Phi_{4,0}^{-1}\Phi_{4,0}(C')$
is birationally equivalent to $V_{14}$, it suffices to see that $\alp$
is generically injective. This follows from the following two facts: first, the
pair $(\xi ,V_{14})$ is determined by $(\G, X)$ uniqueley up to isomorphism,
and second, a generic $V_{14}$ has no biregular automorphisms. If
there were two points $\xi ,\xi '\in V_{14}$ giving the same $\G$, then
there would exist an automorphism of $V_{14}$ sending $\xi$ to $\xi '$,
and hence $\xi =\xi '$. Another proof of the generic injectivity
of $\alp$ is given in Proposition \ref{claim}.

We did not find an appropriate reference for the second fact, so we prove
it in the next lemma.
\end{proof}

\begin{lemma}
A generic variety $V_{14}$ has no nontrivial biregular automorphisms.
\end{lemma}

\begin{proof}
As $V_{14}$ is embedded in $\PP^9$ by the anticanonical system, any
biregular automorphism $g$ of $V_{14}$ is induced by a linear automorphism
of $\PP^9$. Hence it sends conics to conics, and thus defines an
automorphism $F(g):F(V_{14})\lra F(V_{14})$ of the Fano surface $F(V_{14})$,
parametrizing conics on $V_{14}$. In \cite{BD}, the authors prove that
the Hilbert scheme $\mt{Hilb}^2(S)=S^{[2]}$ parametrizing pairs of
points on the K3 surface $S$ of degree 14 in $\PP^8$ is isomorphic
to the Fano 4-fold $F(V_3^4)$ parametrizing lines on $V_3^4$,
where $V_3^4$ is the 4-dimensional linear section of the Pfaffian
cubic in $\PP^{14}$ associated to $S$.
The same argument shows
that $F(V_{14})\simeq F(X)$, where $X$ is the cubic 3-fold associated to $V_{14}$,
and $F(X)$ is the Fano surface parametrizing lines on $X$.

Hence $g$ induces an automorphism $f$ of $F(X)$.
Let $f^*$ be the induced linear automorphism of $\mt{Alb}(F(X))=J^2(X)$,
and $T_0f^*$ its differential at the origin. By \cite{Tyu}, the projectivized tangent
cone of the theta divisor of $J^2(X)$ at $0$ is isomorphic to $X$,
so $T_0f^*$ induces an automorphism of $X$.
$V_{14}$ being generic, $X$ is also generic, so $\mt{Aut} (X)=\{ 1\}$.
Hence $f^*=\id$. By the Tangent Theorem for $F(X)$ \cite{CG},
$\Omega^1_{F(X)}$ is identified with the restriction of the universal
rank 2 quotient bundle $\QQQ$ on $G(2,5)$, and all the global sections of
$\Omega^1_{F(X)}$ are induced by linear forms $L$ on $\PP^4$ via the natural
map $H^0(\PP^4,\OOO_{\PP^4}(1))\otimes\OOO_{G(2,5)}\twoheadrightarrow \QQQ$.
Hence the fact that $f$ acts trivially on $H^0(\Omega^1_{F(X)})=T_0^*J^2(X)$
implies that $f$ permutes the lines $l\in \{ L=0\}\cap X$
lying in one hyperplane section of $X$. For general $L$, there are 27
lines $l$, and in taking two hypeplane sections $\{ L_1=0\} ,\{ L_2=0\}$
which have only one common line, we conclude that $f$ fixes the generic
point of $F(X)$. Hence $F(g)$ is the identity.
This implies that every conic on $V_{14}$ is transformed by $g$
into itself. 

By Theorem \ref{TT}, we have 16 different conics 
$C_1,\ldots ,C_{16}$ passing through the generic point $x\in V_{14}$,
which are transforms of the 16 chords of $C'$ in $X$.
Two different conics $C_i,C_j$ cannot meet at a point $y$, different
from $x$. Indeed, their proper transforms in $X^+$ (we are using
the notations of Theorem \ref{TT}) are the results 
$l_i^+,l_j^+$ of the floppings
of two distinct chords $l_i,l_j$ of $C'$. Two distinct chords
of $C'$ are disjoint, because otherwise the 4 points $(l_i\cap
l_j)\cap C'$ would be coplanar, which would contradict the linear
normality of $C'$. Hence $l_i^+,l_j^+$ are disjoint.
They meet the exceptional divisor $M^+$ of
$X^+\lra V_{14}$ at one point each, hence $C_i\cap C_j=\{ x\}$.
As $g(C_i)=C_i$ and $g(C_j)=C_j$,  this implies $g(x)=x$.
This ends the proof.
\end{proof}

\subsection{Two correspondences between \( H_{4,0},H_{5,1}\)}
This subsection contains some complementary information on the relations
between the families of rational normal quartics and of elliptic
normal quintics on $X$ which can be easily deduced from
the above results.

For a generic cubic 3-fold $X$ and any point $c\in J^2(X)$,
the Abel--Jacobi maps \( \Phi_{4,0},\Phi_{5,1}\)
define a correspondence $Z_1(c)$ between \( H_{4,0},H_{5,1}\)
with generic fibers over $H_{5,1}$, $H_{4,0}$ of dimensions 3,
respectively 5:
$$
Z_1(c)=\{ (\G ,C)\in H_{4,0}\times H_{5,1}\mid
\Phi_{4,0}(\G )+\Phi_{5,1}(C)=c \}\ .
$$
The structure of the fibers is given by Theorems \ref{degree1} and
\ref{fiberquart}: they are, respectiveley, birational to $X$
and isomorphic to $\PP^5$.

There   is another correspondence, defined in \cite{MT}:
\begin{multline*}
Z_2=\{ (\G ,C)\in H_{4,0}\times H_{5,1}\mid
C+\G =\FF_1\cap X\ \mbox{for a rational} \\ \mbox{normal scroll}\
\FF_1\subset\PP^4\}\ .
\end{multline*}
It is proved in \cite{MT} that the fiber over a generic
$C\in H_{5,1}$ is isomorphic to $C$, and the one over $\G\in H_{4,0}$ is
a rational 3-dimensional variety. In fact, we have the following
description for the latter:

\begin{lemma}
For any rational normal quartic $\G\in X$, we have
$Z_2(\G )\simeq PGL(2)$.
\end{lemma}

\begin{proof}
Let $\Gamma\subset\PP^4$ be a rational normal quartic. Then there
exists a unique $PGL(2)$-orbit $PGL(2)\!\cdot\! g\subset PGL(5)$
transforming $\Gamma$ to the normal form
$$
\{ (s^4,s^3t,\ldots ,t^4)\}_{(s:t)\in\PP^1}=
\{ (x_0,\ldots ,x_4)\mid\rk
\left(\begin{array}{cccc}
x_0&x_1&x_2&x_3 \\
x_1&x_2&x_3&x_4
\end{array}\right)
\leq 1\} .
$$
There is one particularly simple rational normal scroll $S$ containing
$\Gamma$:
$$
S=\{ (us^2,ust,ut^2,vs,vt)\}=\{ (x_0,\ldots ,x_4)\mid\rk
\left(\begin{array}{ccc}
x_0&x_1&x_3 \\
x_1&x_2&x_4
\end{array}\right)
\leq 1\} .
$$
Geometrically, $S$ is the union of lines which join the
corresponding points of the line $l=\{ (0,0,0,s,t)\}$
and of the conic $C^2=\{ (s^2,st,t^2$, $0,0)\}$.
Conversely, any rational normal scroll can be obtained
in this way from a pair $(l,C)$ whose linear span is the
whole $\PP^4$. Remark that $(s:t)\mapsto (s:t)$ is the only
correspondence from $l$ to $C$ such that the resulting scroll
contains $\Gamma$.

Now, it is easy to describe all the scrolls containing $\Gamma$:
they are obtained from $S$ by the action of $PGL(2)$. Each non-identical
transformation from $PGL(2)$ leaves invariant $\Gamma$, but
moves both $l$ and $C$, and hence moves $S$.
\end{proof}

As the rational normal scrolls in $\PP^4$ are parametrized by
a rational variety, the Abel--Jacobi image of $C+\G$ is
a constant $c\in J^2(X)$ for all pairs $(\G ,C)$ such that $C\in Z_2(\G )$.
Hence we have identically $\Phi_{4,0}(\G )+\Phi_{5,1}(C)=c$ on $Z_2$,
so that $Z_2(\G )\subset Z_1(c)(\G )$.

We can obtain another birational description of $Z_2(\G )$ for generic $\G$
in applying to all the $C\in Z_2(\G )$ the Tregub--Takeuchi
transformation $\chi$, centered at $\G$. Let $\xi\in V_{14}$ be the
indeterminacy point of $\chi^{-1}$. \smallskip

\begin{proposition}\label{claim}
On a generic $V_{14}$, the family of elliptic 
quintic curves is irreducible. It is parametrized by an open
subset of a component $\BBB$ of $\mt{Hilb}^{5n}_{V_{14}}$ isomorphic to $\PP^5$, and all the curves represented by points of $\BBB$  are l. c. i.
of pure dimension $1$.

For any $x\in V_{14}$, the family of curves from $\BBB =\PP^5$
passing through $x$ is a linear 3-dimensional subspace $\PP^3_x\subset \PP^5$.
For generic rational normal quartic $\G$ as above, $\chi$
maps $Z_2(\G )$ birationally onto $\PP^3_\xi$.
\end{proposition}

\begin{proof}
Gushel constructs in \cite{G2} for any elliptic quintic curve $B$ on $V_{14}$
a rank two vector bundle $\GGG$ such that $h^0 (\GGG )=6$,
$\det \GGG =\OOO (1)$, $c_2 (\GGG )=B$, and proves that the map
from $V_{14}$ to $G=G(2,6)$ given
by the sections of $\GGG$ and composed with the Pl\"ucker embedding
is the standard embedding of $V_{14}$ into $\PP^{14}$.
Hence $\GGG$ is isomorphic to the restriction of  the universal rank
2 quotient bundle on $G$ (in particular, it has no moduli), and the
zero loci of its sections are precisely the sections of $V_{14}$
by the Schubert varieties $\si_{11}(L)$ over all hyperplanes
$L\subset \CC^6=H^0 (\GGG )^\dual$. These zero loci are l. c. i.
of pure dimension 1. Indeed, assume the contrary. Assume
that $D=\si_{11}(L)\cap V_{14}$ has a component of dimension
$>1$. Anyway, $\deg D=\deg \si_{11}(L)=5$, hence if $\dim D=2$,
then $V_{14}$ has a divisor of degree $\leq 5<14=\deg V_{14}$.
This contradicts the fact that $V_{14}$ has index 1 and Picard
number 1. One cannot have $\dim D>2$, because otherwise
$V_{14}$ would be reducible. Hence $\dim D\leq 1$, and it is
l. c. i. of pure dimension 1 as the zero locus of a section
of a rank 2 vector bundle. All the zero loci $B$ of
sections of $\GGG$ form a component $\BBB$ of the Hilbert scheme
of $V_{14}$ isomorphic to $\PP^5$.

The curves $B$ from $\BBB$ passing through
$x$ are the sections of $V_{14}$ by the Schubert varieties
$\si_{11}(L)$ for all $L$ containing the 2-plane $S_x$
represented by the point $x\in G(2,6)$, and hence
form a linear subspace $\PP^3$ in $\PP^5$.

Now, let us prove the last assertion. Let $C\in Z_2(\G )$ be generic.
We have $(C\cdot\G )_{\FF_1}=7$,
therefore the map $\chi$, given by the linear system $|\OOO (8)-5\G )|$,
sends it to a curve $\tilde{C}$ of degree $8\cdot 5-5\cdot 7=5$.
So, the image is a quintic of genus 1. Let $k=\mt{mult}_\xi {\tilde{C}}$.
The inverse $\chi^{-1}$ being given by the linear system
$|\OOO (2)-5\xi |$, we have for the degree of $C=\chi^{-1}({\tilde{C}})$:
\qquad $5=2\deg {\tilde{C}}-5k=10-5k$, hence $k=1$, that is,
$\xi$ is a simple point of ${\tilde{C}}$. Thus the generic $C\in Z_2(\G )$ 
is transformed into a smooth elliptic quintic ${\tilde{C}}\in V_{14}$
passing through $\xi$. By the above, such curves form a $\PP^3$
in the Hilbert scheme, and this ends the proof.
\end{proof}

\subsection{Period map of varieties $V_{14}$}

We have seen that one can associate to any Fano variety $V_{14}$
a unique cubic 3-fold $X$, but to any cubic 3-fold $X$ a 5-dimensional
family of varieties $V_{14}$. Now we are going to determine
this 5-dimensional family. This will give also some
information on the period map of varieties $V_{14}$.
Let ${\mathcal A}_g$ denote the moduli space of 
principally polarized abelian varieties of dimension~$g$.

\begin{theorem}\label{periods}
Let ${\mathcal V}_{14}$ be the moduli space
of smooth Fano $3$-folds of degree $14$, and let
${\Pi}: {\mathcal V}_{14} \rightarrow {\mathcal A}_5$
be the period map on ${\mathcal V}_{14}$.
Then the image ${\Pi}({\mathcal V}_{14})$ coincides with the
$10$-dimensional locus ${\mathcal J}_5$ of intermediate jacobians of cubic
threefolds.
The fiber ${\Pi}^{-1}(J)$, $J \in {\mathcal J}_5$, is isomorphic
to the family ${\mathcal V}(X)$ of the $V_{14}$ which are
associated to the same cubic threefold $X$, and birational
to $J^2(X)$.
\end{theorem}

\begin{proof}
For the construction of ${\mathcal V}_{14}$ and for
the fact that $\dim{\mathcal V}_{14}=15$, see Theorem 0.9 in \cite{Muk}.

According to Theorem \ref{KKK}, there exists a 5-dimensional
family of varieties $V_{14}$, associated to a fixed generic
cubic 3-fold $X$, which
is birationally parametrized by the set $M$ of isomorphism classes of
vector bundle $\EEE$ obtained by Serre's construction
starting from normal elliptic quintics $C\subset X$. By Corollary
\ref{openinJ}, $M$ is an open subset of $J^2(X)$. Hence all the assertions
follow from Proposition \ref{FI-TT}.

\end{proof}


\begin{thebibliography}{-----------}

\bibitem[AC]{A-C}  E. Arrondo, L. Costa, {\em
Vector bundles on Fano 3-folds without intermediate cohomology
}, e-print math.AG/9804033.

\bibitem[AR]{AR} A. Adler, S. Ramanan, {\it Moduli of Abelien Varieties},
Lect. Notes Math., 1644, Springer-Verlag, 1996.

\bibitem[B-MR1]{BM} E. Ballico, R. M. Mir\'o-Roig,
{\it Rank $2$ stable vector bundles on Fano $3$-folds of index $2$},
J. Pure Appl. Algebra {\bf 120}, 213--220 (1997).

\bibitem[B-MR2]{BM2} E. Ballico, R. M. Mir\'o-Roig,
{\it A lower bound for the number of components
of the moduli schemes of stable rank $2$ vector bundles on projective $3$-folds. }, Proc. Amer.
Math. Soc. {\bf 127}, 2557--2560 (1999).

\bibitem[Barth]{Barth} W. Barth, {\it 
Irreducibility of the space of mathematical instanton bundles with rank $2$ and
$c\sb{2}=4$},
Math. Ann. {\bf 258}, 81--106 (1981/82).  



\bibitem[BS]{BS} D. Bayer, M. Stillman: {\em Macaulay: A computer
algebra system for algebraic geometry}, 1982+, ftp://zariski.harvard.edu
(login: anonymous Password: any).

\bibitem[B]{B} A. Beauville,
{\it Les singularit\'es du diviseur th\^eta
de la jacobienne interm\'ediaire de l'hypersurface
cubique dans ${\bf P}^{4}$,}
Lecture Notes in Math. {\bf 947}, 190-208 (1982).

\bibitem[BD]{BD} A. Beauville,
R. Donagi:
{\it La vari\'et\'e des droites d'une hypersurface cubique de dimension~$4$}, 
C. R. Acad. Sci. Paris Ser. I Math. {\bf  301}, no. 14, 703--706(1985).


\bibitem[CG]{CG} H.Clemens, P.Griffiths,
{\it The intermediate
Jacobian of the cubic threefold},
Annals of Math. {\bf 95}, 281-356 (1972).



\bibitem[ES]{ES} G. Ellingsrud, S. Str\o mme, {\it 
Stable rank-$2$ vector bundles on ${P}\sp{3}$ with
$c\sb{1}=0$ and $c\sb{2}=3$}, 
Math. Ann. {\bf 255}, 123--135 (1981).

\bibitem[E-SB]{E-SB} L.Ein, N. Shepherd-Barron, {\em Some special Cremona transformations}, Amer. J. Math. {\bf 111}, 783-800 (1989).



\bibitem[F]{F} G. Fano, {\it Sulle sezioni spaziale della variet\`a
grassmanniana delle rette spazio a cinque dimensioni},
Rend. R. Accad. Lincei {\bf 11}, 329--335(1930).

\bibitem[G1]{G1} N. Gushel, {\it Fano varieties of genus $6$},
Math. USSR, Izv. {\bf 21}, 445--459
(1983).

\bibitem[G2]{G2} N. Gushel, {\it Fano 3-folds of genus $8$},
St. Petersbg. Math. J.  {\bf 4}, 115--129
(1993).


\bibitem[H]{H} R.Hartshorne,
{\it 
Stable vector bundles of rank $2$ on ${P}\sp{3}$}, 
Math. Ann. {\bf 238}, 229--280 (1978). 



\bibitem[HH]{HH} Hartshorne, R.; Hirschowitz, A. 
{\em Smoothing algebraic space 
curves,} In: Algebraic geometry, Sitges (Barcelona), 1983, 
Lecture Notes in Math., 1124. Springer, Berlin-New
York, 1985,  98--131.


\bibitem[Hu]{Hu} Hulek, K.: {\em Projective geometry of elliptic
curves,} Ast\'erisque {\bf 137}(1986).


\bibitem[I]{Ili} A.I.Iliev,
{\it Minimal sections of conic bundles},
Boll. Uni. Mat. Italiana. (8) {\bf 2B}, 401-428(1999).

\bibitem[IR]{IR} A.I.Iliev, K.Ranestad,
{\it The sympletic grassmannian $Sp(3)/U(3)$
and duality for the Fano threefold $V_{16}$},
e-print math.AG/9809140, to appear in J. Alg. Geom.

\bibitem[Isk1]{Isk-1} V.A.Iskovskikh,
{\it Birational automorphisms of
three-dimensional algebraic varieties},
J. Soviet Math. {\bf 13:6}, 815-868 (1980)

\bibitem[Isk2]{Isk-2} V.A. Iskovskikh,
{\it Double projection from a line on
Fano threefolds of the first kind},
Math. USSR Sbornik {\bf 66}, No.1,
265-284 (1990).

\bibitem[Isk-P]{IP} V. A. Iskovskikh, Yu. G. Prokhorov,
{\em Fano varieties}, Algebraic Geometry V (A. N. Parsin,
I. R .Shafarevich, Eds.), Encyclopaedia of Math. Sci., Vol. 47,
Springer, Berlin, 1998.




\bibitem[LP]{LP} J. Le Potier, {\it
Sur l'espace de modules des fibr\'es de Yang et Mills}, In:
Mathematics and physics
(Paris, 1979/1982), 65--137,  
Progr. Math., 37, Birkh\"auser, Boston, Mass., 1983.

\bibitem[M]{M} D. Markushevich, {\it  Numerical invariants of families of 
straight lines of Fano manifolds}, Math.USSR, 
Sbornik {\bf 44}, 239-260(1983).

\bibitem[MT]{MT} D. Markushevich, A.S.Tikhomirov,
{\it The Abel-Jacobi map of a moduli component
of vector bundles on the cubic threefold}, e-print math.AG/9809140,
to appear in J. Alg. Geom.


\bibitem[Mu]{Muk} S.Mukai,
{\it Curves, K3 Surfaces and Fano 3-folds of Genus $\le$ 10},
in:
{\it Algebraic Geometry and Commutative Algebra in Honor of
M. Nagata}, Kinokuniya, Tokyo 357-377 (1988)

\bibitem[OS]{OSz} G. Ottaviani, M. Szurek,
{\em On moduli of stable $2$-bundles with small Chern classes
on $Q_3$}, Ann. Mat. Pura Appl. (IV) {\bf 167}, 191--241 (1994).

\bibitem[P]{Pu} P.J. Puts,
{\it On some Fano threefolds that are
sections of Grassmannians},
Proc. Kon. Ned. Acad. Wet. {\bf A 85(1)},
77-90 (1982)



\bibitem[SW]{SW} M. Szurek, J. Wi\'sniewski, {\em
Conics, conic fibrations and
stable vector bundles of rank $2$ on some Fano threefolds}, 
Rev. Roumaine Math. Pures Appl. {\bf 38}, 729--741 (1993). 

\bibitem[Tak]{Tak} K.Takeuchi,
{\it Some birational maps of Fano $3$-folds},
Compositio Math. {\bf 71}, 265-283 (1989).

\bibitem[Tre]{Tre} S.L.Tregub,
{\it Construction of a birational isomorphism of a cubic
threefold and Fano variety of the first kind with $g=8$, associated
with a normal rational curve of degree 4},
Moscow Univ. Math. Bull.
{\bf 40},  78-80 (1985); translation from Vestn. Mosk. Univ., Ser. I
1985, No.6, 99-101 (1985).


\bibitem[Tyu]{Tyu} A. N. Tyurin, {\em Geometry 
of the Fano surface of a non-singular cubic 
$F\subset \PP^4$ and Torelli theorems for Fano surfaces 
and cubics,} Izv. USSR Math. {\bf 35} (1971), 498--529.


\end{thebibliography}
\end{document}